\input amstex
\magnification=\magstep1 
\baselineskip=13pt
\documentstyle{amsppt}
\vsize=8.7truein \CenteredTagsOnSplits \NoRunningHeads
 \def\EE{\bold{E}\thinspace}
 \def\PP{\bold{P}\thinspace}
 \def\tr{\operatorname{trace}}
 \def\rk{\operatorname{rank}}
 \def\conv{\operatorname{conv}}
 
 \def\sym{\operatorname{Sym}}
 \topmatter
 
 \title Convexity of the image of a quadratic map via the relative entropy distance\endtitle
 \author Alexander Barvinok \endauthor
 \address Department of Mathematics, University of Michigan, Ann Arbor,
MI 48109-1043, USA \endaddress
\email barvinok$\@$umich.edu \endemail
\date May 2013 \enddate
 \thanks  This research was partially supported by NSF Grant DMS 0856640.
\endthanks 
 \abstract Let $\psi: {\Bbb R}^n \longrightarrow {\Bbb R}^k$ be a map defined by $k$ positive definite quadratic forms on 
 ${\Bbb R}^n$. We prove that the relative entropy (Kullback-Leibler) distance from the convex hull of the image of $\psi$ to the image of 
 $\psi$ is bounded above by an absolute constant. More precisely, we prove that for every point $a=\left(a_1, \ldots, a_k\right)$ in the convex hull of the image of $\psi$ 
 such that $a_1 + \ldots + a_k=1$ there is a point $b=\left(b_1, \ldots, b_k\right)$ in the image of $\psi$ such that 
 $b_1 + \ldots + b_k =1$ and such that $\sum_{i=1}^k a_i \ln \left(a_i/b_i\right) < 4.8$. Similarly, we prove that for any integer
 $m$ one can choose a convex combination $b$ of at most $m$ points from the image of $\psi$ such that 
 $\sum_{i=1}^k a_i \ln \left(a_i/b_i\right) < 15/\sqrt{m}$.
 \endabstract

\endtopmatter

\document

\head 1. Introduction \endhead

Let $q_1, \ldots, q_k: {\Bbb R}^n \longrightarrow {\Bbb R}$ be quadratic forms and let 
$\psi: {\Bbb R}^n \longrightarrow {\Bbb R}^k$ be the corresponding quadratic map,
$$\psi(x)=\left(q_1(x), \ldots, q_k(x)\right).$$
We are interested in the convex properties of the image $\psi\left({\Bbb R}^n\right) \subset {\Bbb R}^k$. The image is clearly convex when $k=1$ and by the Dines Theorem it is convex when $k=2$ (this and related facts can be found, for example, in 
Sections II.12-14 of \cite{Ba02} or in \cite{PT07}). The image is not necessarily convex for $k \geq 3$, though it remains convex for $k=3$ if some linear combination of the forms $q_1, q_2$ and $q_3$ is positive definite.

In this paper, we show that the image $\psi\left({\Bbb R}^n\right)$ is close to its own convex hull 
$\conv\left( \psi\left({\Bbb R}^n\right)\right)$ in some information-theoretic sense. 

Let 
$$a=\left(a_1, \ldots, a_k\right) \quad \text{and} \quad b=\left(b_1, \ldots, b_k \right)$$
be two positive vectors such that 
$$\sum_{i=1}^k a_i =\sum_{i=1}^k b_i =1.$$
We interpret $a$ and $b$ as probability distributions and define the {\it relative entropy} of $a$ with respect to $b$ as 
$$D(a\|b)=\sum_{i=1}^k a_i \ln \left({ a_i \over b_i}\right).$$
The quantity $D(a\|b)$ is also known as the {\it Kullback - Leibler distance} from $a$ to $b$ (although, generally speaking,
$D(a\|b)\ne D(b\|a)$ and the triangle inequality does not hold).
In particular, $D(a\|b) \geq 0$ with equality if and only if $a=b$, see for example, \cite{CT06}. 

We prove that with respect to the Kullback - Leibler distance, the image $\psi\left({\Bbb R}^n\right)$ of a quadratic map is 
reasonably close to its own convex hull $\conv\left(\psi\left({\Bbb R}^n\right)\right)$.
\proclaim{(1.1) Theorem} Let $q_1, \ldots, q_k: {\Bbb R}^n \longrightarrow {\Bbb R}$ be positive definite quadratic forms and 
let $\psi: {\Bbb R}^n \longrightarrow {\Bbb R}^k$ be the corresponding map,
$$\psi(x)=\left(q_1(x), \ldots, q_k(x) \right).$$
Let $a \in \conv\left( \psi\left({\Bbb R}^n \right)\right)$ be a point, $a=\left(a_1, \ldots, a_k\right)$, such that 
$a_1 + \ldots +a_k=1$.
Then there exists a point $b \in \psi\left({\Bbb R}^n \right)$, $b=\left(b_1, \ldots, b_k\right)$, such that 
$b_1 +\ldots +b_k=1$ and 
$$\sum_{i=1}^k a_i \ln \left({a_i \over b_i} \right) \ \leq \ \beta$$
for some absolute constant $\beta>0$. One can choose, for example, $\beta =4.8$.
\endproclaim
We have undertaken some effort to optimize the constant $\beta$, but its optimal value is not known at the moment and it would be interesting to find it.

Loosely speaking, Theorem 1.1 asserts that replacing the image of $\psi$ by its convex hull leads to only a constant loss of information. The technique of {\it semidefinite programming} is based on replacing computationally intractable systems of quadratic equations and inequalities over the reals by computationally tractable systems of linear equations and inequalities in positive semidefinite matrices. This procedure is known as {\it relaxation}, see for example, \cite{Tu10}. The success of relaxation depends on the convex properties of the underlying quadratic maps, see
\cite{PT07}. Speaking even more loosely, one can speculate that the constant bound on the information loss in Theorem 1.1 explains the success of semidefinite programming.

We also prove the following extension of Theorem 1.1.

\proclaim{(1.2) Theorem} Let $q_1, \ldots, q_k: {\Bbb R}^n \longrightarrow {\Bbb R}$ be positive definite quadratic forms and 
let $\psi: {\Bbb R}^n \longrightarrow {\Bbb R}^k$ be the corresponding map,
$$\psi(x)=\left(q_1(x), \ldots, q_k(x) \right).$$
Let $a \in \conv\left( \psi\left({\Bbb R}^n \right)\right)$ be a point, $a=\left(a_1, \ldots, a_k\right)$, such that 
$a_1 + \ldots +a_k=1$.
Then, for any positive integer $m$, there exists a point $b=\left(b_1, \ldots, b_k\right)$, such that 
$b_1 +\ldots +b_k=1$, the point $b$ is a convex combination of at most $m$ points of $\psi\left({\Bbb R}^n\right)$ and 
$$\sum_{i=1}^k a_i \ln \left({a_i \over b_i} \right) \ < \ {15 \over \sqrt{m}}.$$
\endproclaim

We note a useful inequality
$$D(a\|b) =\sum_{i=1}^k a_i \ln \left({a_i \over b_i}\right) \ \geq \ {1 \over 2 \ln 2} \left(\sum_{i=1}^k \left| a_i -b_i \right| \right)^2,$$
see, for example, Section 11.6 of \cite{CT06}.
The {\it Approximate Carath\'eodory Theorem} of Maurey (see \cite{Pi81} and Section I.3 of \cite{Ve+}) states that if $X$ is 
{\it any} set of points in the standard simplex 
$$\sum_{i=1}^k x_i =1 \quad \text{and} \quad x_1, \ldots, x_k \geq 0$$ 
in ${\Bbb R}^k$ then any point $a \in \conv(X)$ can be approximated within error of $1/\sqrt{m}$ by a convex combination of $m$ points of $X$ in the $\ell^2$ (Euclidean) norm. Theorem 1.2 asserts that if $X$ is the image of a quadratic map then one can get a similar approximation in the $\ell^1$ norm.

The {\it Johnson - Lindenstrauss Lemma} implies that for any $\epsilon>0$, if one chooses $m=O\left(\epsilon^{-2} \ln k \right)$ in Theorem 1.2 then one can ensure that
$$\left| \ln {a_i \over b_i} \right| \ \leq \ \epsilon \quad \text{for} \quad i=1, \ldots, k,$$
see, for example, Sections V.5-6 of \cite{Ba02} and \cite{Ma08}. Theorem 1.2 asserts that if we measure the Kullback - Leibler distance, then 
the dependence on the number $k$ of quadratic forms can be removed so that $m=O\left(\epsilon^{-2}\right)$ and
$$D(a\|b)=\sum_{i=1}^k a_i \ln \left({a_i \over b_i}\right) \ \leq \ \epsilon.$$

In the rest of the paper, we prove Theorems 1.1 and 1.2. In Section 2, we establish some general results on the distribution of values of a positive semidefinite quadratic form with respect to the Gaussian probability measure in ${\Bbb R}^n$. In Section 3, we consider 
the problem of maximizing a convex combination of logarithms of positive semidefinite quadratic forms on the unit sphere. We prove that its straightforward positive semidefinite relaxation produces a relative error bounded by an absolute constant. In Section 4, we complete the proof of Theorem 1.1. The proof of Theorem 1.2 given in Section 5 is a straightforward modification of our proof of Theorem 1.1.

\head 2. Quadratic forms and the Gaussian measure \endhead

Un this section, we prove the following main result.
\proclaim{(2.1) Lemma} Let us fix in ${\Bbb R}^n$ the standard Gaussian probability measure $\mu_n$ with density 
$${1 \over (2\pi)^{n/2}} e^{-\|x\|^2/2}.$$
Let $q: {\Bbb R}^n \longrightarrow {\Bbb R}$ be a positive semidefinite quadratic form such that 
$$\EE q=1.$$
Then
\roster
\item We have $$\EE \left| \ln q  \right| \ < \ 2.75;$$
\item For $t \geq 1$ let us define
$$\phi(t) =\min_{\alpha \geq 1} {2^{\alpha} \over t^{\alpha} \sqrt{\pi}} \Gamma\left(\alpha +{1 \over 2}\right).$$
Then
$$\PP\bigl(x:\ q(x) \geq t \bigr) \ \leq \ \phi(t) \quad \text{for all} \quad t \geq 1.$$
\endroster
\endproclaim
\demo{Proof} Part (1) is essentially proved in \cite{Ba99} but we present its proof here for completeness. 
We have 
$$\EE \left| \ln q \right| \ \leq \ \left(\EE \ln^2 q \right)^{1/2}.$$
We can write
$$q(x)=\sum_{i=1}^n \lambda_i x_i^2 \quad \text{for} \quad x=\left(x_1, \ldots, x_n \right) \tag2.1.1$$
in some orthonormal basis of ${\Bbb R}^n$. Since 
$$\EE q =\EE x_i^2 =1 \quad \text{for} \quad i=1, \ldots, n,$$
we have 
$$\sum_{i=1}^n \lambda_i=1 \qquad \text{and also} \quad \lambda_i \geq 0 \quad \text{for} \quad i=1, \ldots, n. \tag2.1.2$$
Let
$$Y=\Bigl\{x \in {\Bbb R}^n: \ q(x) \leq 1 \Bigr\}.$$
By the concavity of the logarithm,
$$\ln \left( \sum_{i=1}^n \lambda_i x_i^2 \right) \ \geq \ \sum_{i=1}^n \lambda_i \ln x_i^2.$$
Since $\ln q(x) < 0$ for all $x \in Y$, using (2.1.2) and the convexity of the function $t \longmapsto t^2$, we conclude that
$$\split \int_Y \ln^2 q(x) \ d\mu_n(x) \ \leq \ &\int_Y \left( \sum_{i=1}^n \lambda_i  \ln x_i^2 \right)^2 \ d\mu_n(x) 
\\ \leq \ &\int_Y \left(\sum_{i=1}^n \lambda_i \ln^2 x_i^2 \right) \ d \mu_n(x)
\leq  \int_{{\Bbb R}^n} \ln^2 x_1^2 \ d\mu_n(x) \\ =&{8 \over \sqrt{2 \pi}} \int_0^{+\infty}\left( \ln^2 x \right) e^{-x^2/2} \ dx  \ < \ 6.55. \endsplit$$
On the other hand, since $\ln t \ \leq \ \sqrt{t}$ for $t \geq 1$, we conclude that 
$$\int_{{\Bbb R}^n \setminus Y} \ln^2 q(x) \ d \mu_n(x) \ \leq \ \int_{{\Bbb R}^n\setminus Y} q(x) \ d\mu_n(x)\  \leq \ \int_{{\Bbb R}^n} q(x) \ d \mu_n(x) =1.$$
Therefore,
$$\EE \ln^2 q \ < \ 6.55 +1 =7.55 \quad \text{and} \quad \EE \left| \ln q \right| \ < \ \sqrt{7.55} \ < \ 2.75,$$
which proves Part (1).

Let us choose any $\alpha \geq 1$. Applying the Markov inequality, we get
$$\PP\bigl(x:\ q(x) \geq t\bigr)=\PP\left(x: \ q^{\alpha}(x) \geq t^{\alpha} \right) \ \leq \ t^{-\alpha} \EE q^{\alpha}.$$
Writing $q$ as in (2.1.1) and using (2.1.2) and the convexity of the function $t \longmapsto t^{\alpha}$, we obtain
$$\split \EE q^{\alpha} =&\EE\left( \sum_{i=1}^n \lambda_i x_i^2 \right)^{\alpha} \ \leq \ \sum_{i=1}^k \lambda_i \EE \left(x_i^{2}\right)^{ \alpha}\\= &{2 \over \sqrt{2 \pi}} \int_0^{+\infty} x^{2\alpha} e^{-x^2/2} \ dx = 
{2^{\alpha} \over \sqrt{\pi}} \Gamma\left(\alpha +{1 \over 2}\right), \endsplit$$
from which the proof of Part (2) follows.
{\hfill \hfill \hfill} \qed
\enddemo

\remark{(2.2) Remark} The exact upper bound in Part (1) is not known to the author, though it looks plausible that it is  attained on forms of rank 1 and hence is equal to
$${4 \over \sqrt{2 \pi}} \int_0^{+\infty}|\ln x| e^{-x^2/2} \ dx \approx 1.76.$$
\endremark

\head 3. An optimization problem on the sphere \endhead 

\subhead (3.1) Notation \endsubhead
We consider the space $\sym_n$ of $n \times n$ symmetric matrices endowed with standard inner product 
$$\langle A, B \rangle =\sum_{i,j=1}^n a_{ij} b_{ij} = \tr(AB),$$
where $A=\left(a_{ij}\right)$ and $B=\left(b_{ij}\right)$. For a vector $x \in {\Bbb R}^n$,
$x=\left(x_1, \ldots, x_n \right)$, we define a symmetric matrix 
$X =x \otimes x$, $X=\left(x_{ij}\right)$, by $x_{ij}=x_i x_j$. Thus a quadratic form $q$ with matrix $Q$ can be written 
as 
$$q(x)=\langle Q,\ x \otimes x \rangle \quad \text{for all} \quad x \in {\Bbb R}^n.$$ We write $X \succeq 0$ to say that $X$ is 
positive semidefinite and $X \succ 0$ to say that $X$ is positive definite.

In ${\Bbb R}^n$, we consider the standard inner product 
$$\langle x, y \rangle =\sum_{i=1}^n x_i y_i \quad \text{where} \quad x=\left(x_1, \ldots, x_n\right) \quad \text{and} \quad 
y=\left(y_1, \ldots, y_n \right),$$
the corresponding norm 
$$\|x\|=\sqrt{\langle x, x \rangle},$$
and the unit sphere
$${\Bbb S}^{n-1}=\Bigl\{x \in {\Bbb R}^n:\quad \|x\| =1 \Bigr\}.$$

In this section, we prove the following main result.

\proclaim{(3.2) Theorem} Let $\alpha_1, \ldots, \alpha_k$ be non-negative reals such that 
$\alpha_1+ \ldots + \alpha_k=1$, let $Q_1, \ldots, Q_k$ be $n \times n$ positive definite matrices and let 
$q_1, \ldots, q_k: {\Bbb R}^n \longrightarrow {\Bbb R}$ be the corresponding quadratic forms,
$$q_i(x) =\langle Q_i,\ x \otimes x \rangle \quad \text{for} \quad i=1, \ldots, k.$$
Then
$$\max_{x \in {\Bbb S}^{n-1}} \sum_{i=1}^k \alpha_i \ln q_i(x) \ \leq \ 
\max \Sb X \succeq 0 \\ \tr(X)=1 \endSb \sum_{i=1}^k \alpha_i \ln \langle Q_i, X \rangle \ \leq \ 
\beta + \max_{x \in {\Bbb S}^{n-1}} \sum_{i=1}^k \alpha_i \ln q_i(x),$$
where $\beta>0$ is an absolute constant. One can choose $\beta=4.8.$
\endproclaim
\demo{Proof} For $x \in {\Bbb S}^{n-1}$ the matrix $X =x \otimes x$ satisfies the constraints $X \succeq 0$ and $\tr(X)=1$. Hence the first inequality holds.

Let $A$ be a matrix where the maximum value of the function 
$$X \longmapsto \sum_{i=1}^k \alpha_i \ln \langle Q_i, X \rangle$$
is attained on the set $X$ of positive semidefinite matrices of trace 1.
Rescaling $Q_i \longrightarrow \tau_i Q_i$ for some positive $\tau_1, \ldots, \tau_k$ if necessary, we may assume 
that $\langle Q_i , A \rangle =1$ for $i=1, \ldots, k$ and hence
$$\max \Sb X \succeq 0 \\ \tr(X)=1 \endSb \sum_{i=1}^k \alpha_i \ln \langle Q_i, X \rangle =
 \sum_{i=1}^k \alpha_i \ln \langle Q_i, A \rangle =0. \tag3.2.1$$
 Since $A$ is positive semidefinite, we can write $A=T^2$ for some symmetric
$n \times n$ matrix $T$. 

Let us fix the standard Gaussian probability measure $\mu_n$ in ${\Bbb R}^n$ with density 
$${1 \over (2\pi)^{n/2}} e^{-\|x\|^2/2}$$
and let $x \in {\Bbb R}^n$ be a random vector.
Then 
$$\EE \| Tx\|^2 =\EE \langle Tx,\ Tx \rangle = \EE \langle T^2 x,\ x \rangle = \tr\left(T^2 \right)=
\tr(A)=1.$$
Hence by Part (2) of Lemma 2.1,
$$\PP\bigl(x:\ \|Tx\|^2 \geq \ 6 \bigr)\ \leq \ \phi(6) \ < \ 0.07 \tag3.2.2$$
(choosing $\alpha=3$ in the definition of $\phi(6)$, we obtain $\phi(6) \leq 5/72 < 0.07$).

Furthermore,
$$\split \EE q_i(Tx)= &\langle Q_i Tx,\ Tx \rangle =\langle T Q_i T x,\ x \rangle =\tr \left(T Q_i T \right)\\=
&\tr \left(Q_i T^2 \right) =\langle Q_i, A \rangle =1 \quad \text{for} \quad i=1, \ldots, k. \endsplit$$
Therefore, by Part (1) of Lemma 2.1, 
$$\EE \left| \ln q_i(Tx)\right| \ \leq \ 2.75 \quad \text{for} \quad i=1, \ldots, k$$ 
and hence 
$$\EE \left| \sum_{i=1}^k \alpha_i \ln q_i(Tx) \right| \ \leq \ 2.75.$$
Therefore, by the Markov inequality,
$$\PP \left(x:\quad \sum_{i=1}^k \alpha_i \ln q_i(Tx) \ \leq \ -3 \right) \ \leq \ {2.75 \over 3} \ < \ 0.92. \tag3.2.3$$
From (3.2.2)--(3.2.3) we conclude that there is an $x \in {\Bbb R}^n \setminus \{0\}$ such that 
$$\|Tx\|^2 \ < \ 6 \quad \text{and} \quad \sum_{i=1}^k \alpha_i \ln q_i(Tx) \ > \ -3.$$
Then for 
$$y= {Tx \over \|Tx\|}$$
we have
$$y \in {\Bbb S}^{n-1} \quad \text{and} \quad \sum_{i=1}^k \alpha_i \ln q_i(y) \ > \ -3-\ln(6) \ > -4.8,$$
and, in view of (3.2.1), the proof follows.
{\hfill \hfill \hfill} \qed
\enddemo

\head 4. Proof of Theorem 1.1 \endhead

\demo{Proof} Let us write
$$q_i(x) =\langle Q_i,\ x \otimes x \rangle \quad \text{for} \quad i=1, \ldots, k,$$
where $Q_1, \ldots, Q_k$ are $n \times n$ positive definite matrices. Let 
$$S=\sum_{i=1}^k Q_i.$$
Thus $S \succ 0$ and hence there exists an invertible symmetric matrix $T: {\Bbb R}^n \longrightarrow {\Bbb R}^n$ 
such that $S=T^2$. Let us define new matrices 
$$\widehat{Q}_i = T^{-1} Q_i T^{-1} \quad \text{for} \quad i=1, \ldots, k,$$
the corresponding quadratic forms 
$$\widehat{q}_i(x) = \langle \widehat{Q}_i,\  x \otimes x \rangle = \langle Q_i,\ T^{-1} x \otimes T^{-1} x \rangle =
q_i\left(T^{-1}x\right)
\quad \text{for} \quad i=1, \ldots, k$$
and the map $\widehat{\psi}: {\Bbb R}^n \longrightarrow {\Bbb R}^k$,
$$\widehat{\psi}(x)=\left(\widehat{q}_1, \ldots, \widehat{q}_k \right).$$
Clearly, $\psi\left({\Bbb R}^n \right) =\widehat{\psi}\left({\Bbb R}^n\right)$ and 
$$\sum_{i=1}^k \widehat{Q}_i =I.$$
Hence, without loss of generality, we can assume that 
$$\sum_{i=1}^k Q_i =I. \tag4.1$$
Since $a \in \conv\left(\psi\left({\Bbb  R}^n\right)\right)$, we can write 
$$a_i =\langle Q_i, X \rangle \quad \text{for} \quad i=1, \ldots, k$$ and 
some $X \succeq 0$. Moreover, in view of (4.1), we have 
$$1=\sum_{i=1}^k a_i = \left\langle \sum_{i=1}^k Q_i,\ X \right\rangle = \langle I,\ X \rangle =\tr(X).$$
We note that 
$$\sum_{i=1}^k a_i \ln  \langle Q_i, X \rangle = \sum_{i=1}^k a_i \ln a_i.$$
By Theorem 3.2, there is an $x \in {\Bbb S}^{n-1}$ such that 
$$\beta + \sum_{i=1}^k a_i \ln q_i(x) \ \geq \  \sum_{i=1}^k a_i \ln a_i.$$
Letting 
$$b_i = q_i(x) \quad \text{for} \quad i=1, \ldots, k,$$
we conclude that 
$$\sum_{i=1}^k b_i = \sum_{i=1}^k \langle Q_i,\ x \otimes x \rangle = \langle I, x\otimes x \rangle = 
\tr(x \otimes x) = 1$$ 
and that 
$$\sum_{i=1}^k a_i \ln \left({a_i \over b_i} \right) =\sum_{i=1}^k a_i \ln a_i - \sum_{i=1}^k a_i \ln b_i  \ \leq \ \beta.$$
Moreover, for $b=\left(b_1, \ldots, b_k\right)$ we have $b=\psi(x)$, so $b \in \psi\left({\Bbb R}^n\right)$.
{\hfill \hfill \hfill} \qed
\enddemo

\head 5. Proof of Theorem 1.2 \endhead

\proclaim{(5.1) Lemma} For a positive integer $m$ let us consider ${\Bbb R}^{mn}$ as the direct sum 
$${\Bbb R}^{mn} = \underbrace{{\Bbb R}^n \oplus \ldots \oplus {\Bbb R}^n}_{\text{$m$ times}}.$$
Let us fix the standard Gaussian probability measure $\mu_n$ in ${\Bbb R}^n$ and consider the standard Gaussian probability 
measure $\mu_{mn}$ in ${\Bbb R}^{mn}$ as the direct product
$$\mu_{mn} =\mu_n \otimes \ldots \otimes \mu_n.$$
Let $q: {\Bbb R}^n \longrightarrow {\Bbb R}$ be a positive semidefinite quadratic form and let us define a quadratic form
$q_m: {\Bbb R}^{mn} \longrightarrow {\Bbb R}$ by 
$$q_m\left(x_1, \ldots, x_m \right) ={1 \over m} 
\sum_{i=1}^m q\left(x_i\right) \quad \text{where} \quad x=\left(x_1, \ldots, x_m\right)$$
and $x_i \in {\Bbb R}^n$ for $i=1, \ldots, m$.
Suppose that 
$$\EE q=1.$$
Then
\roster
\item For all $t \geq 1$ we have 
$$\PP\Bigl(x \in {\Bbb R}^{mn}:\ q_m(x) \geq t \Bigr) \ \leq \ \exp\left\{{m \over 2} \left(1-t + \ln t\right)\right\};$$
\item For all $0 < t \leq 1$ we have 
$$\PP\Bigl(x \in {\Bbb R}^{mn}:\ q_m(x) \leq t\Bigr) \ \leq \  \exp\left\{ {m \over 2} \left(1-t+\ln t \right) \right\};$$
\item We have
$$\EE \left| \ln q_m \right| \ \leq \ {6 \over \sqrt{m}}.$$
\endroster
\endproclaim
\demo{Proof} We use the Laplace transform method, see also \cite{HW71}.
Since 
$$\EE q=1,$$
in some orthonormal basis of ${\Bbb R}^{n}$ we can write 
$$q(x)= \sum_{i=1}^n \lambda_i \xi_i^2 \quad \text{where} \quad x=\left(\xi_1, \ldots, \xi_n\right)$$ 
and 
$$\sum_{i=1}^n \lambda_i =1 \quad \text{and} \quad \lambda_i \geq 0 \quad \text{for} \quad i=1, \ldots, n. \tag5.1.1$$
Writing vectors $x \in {\Bbb R}^{mn}$ as $x=\left(\xi_{11}, \ldots, \xi_{1n}, \xi_{21}, \ldots, \xi_{2n}, \ldots, 
\xi_{m1}, \ldots, \xi_{mn} \right)$,
we write 
$$q_m(x)={1 \over m} \sum_{i=1}^n \sum_{j=1}^m \lambda_i \xi_{ji}^2.$$
For any $0 < \alpha < m/2$ we have 
$$\split \PP\Bigl(x \in {\Bbb R}^{mn}: \ q_m(x) \geq t \Bigr) =&\PP\Bigl(x \in {\Bbb R}^{mn}: \ e^{\alpha q_m(x)} \geq e^{\alpha t }\Bigr)
\ \leq \ e^{-\alpha t} \EE e^{\alpha q_m}\\
=&e^{-\alpha t} \prod_{i=1}^n \left( 1 -{2 \alpha \lambda_i \over m}  \right)^{-m/2}. \endsplit$$
Since the function 
$$\left(\lambda_1, \ldots, \lambda_n\right) \longmapsto -{m \over 2} \sum_{i=1}^n \ln \left(1 -{2 \alpha \lambda_i \over m}\right)$$
is convex, it attains its maximum on the simplex (5.1.1) at a vertex $\lambda_i=1$, $\lambda_j=0$ for $j \ne i$. 
Therefore,
$$\PP\Bigl(x \in {\Bbb R}^{mn}: \ q_m(x) \geq t \Bigr)  \ \leq \ e^{-\alpha t} \left(1- {2 \alpha \over m}\right)^{-m/2}.$$
Optimizing on $\alpha$, we choose 
$$\alpha={m \over 2} \left({t-1 \over t}\right)$$
and the proof of Part (1) follows.

For any $\alpha >0$ we have 
$$\split \PP\Bigl(x \in {\Bbb R}^{mn}: \ q_m(x) \leq t \Bigr) =&\PP\Bigl(x \in {\Bbb R}^{mn}: \ e^{-\alpha q_m(x)} \geq e^{-\alpha t }\Bigr)
\ \leq \ e^{\alpha t} \EE e^{-\alpha q_m}\\
=&e^{\alpha t} \prod_{i=1}^n \left( 1 +{2 \alpha \lambda_i \over m}  \right)^{-m/2}. \endsplit$$
Since the function 
$$\left(\lambda_1, \ldots, \lambda_n\right) \longmapsto -{m \over 2} \sum_{i=1}^n \ln \left(1 +{2 \alpha \lambda_i \over m}\right)$$
is convex, it attains its maximum on the simplex (5.1.1) at a vertex $\lambda_i=1$, $\lambda_j=0$ for $j \ne i$. 
Therefore,
$$\PP\Bigl(x \in {\Bbb R}^{mn}: \ q_m(x) \leq t \Bigr)  \ \leq \ e^{\alpha t} \left(1+ {2 \alpha \over m}\right)^{-m/2}.$$
Optimizing on $\alpha$, we choose 
$$\alpha={m \over 2} \left({1-t \over t}\right)$$
and the proof of Part (2) follows.

Let us define 
$$X_+=\Bigl\{ x \in {\Bbb R}^{mn}:\ q_m(x) \geq 1 \Bigr\}\quad \text{and} 
\quad  X_-=\Bigl\{ x \in {\Bbb R}^{mn}:\ q_m(x) < 1 \Bigr\} $$
Then
$$ \EE \left| \ln q_m \right| = \int_{X_+}  \ln q_m(x) \ d\mu_{mn}(x) - \int_{X_-}  \ln q_m(x)\ d\mu_{mn}(x)$$
By Part (1),
$$\split \int_{X_+}  \ln q_m(x) \ d\mu_{mn}(x)=&\int_0^{+\infty} \PP \Bigl(x:\ \ln q_m(x) \geq t \Bigr) \ dt \\=&\int_0^{+\infty} \PP \Bigl(x:\ q_m(x) \geq e^t \Bigr) \ dt \\ \leq \ &\int_0^{+\infty} \exp\left\{ {m \over 2} 
\left(1- e^t + t \right) \right\} \ dt \ \leq \ \int_0^{+\infty} \exp\left\{ -{m t^2 \over 4} \right\} dt\\ =&\sqrt{\pi \over m}. \endsplit$$
By Part(2),
$$\split \int_{X_-}  -\ln q_m(x) \ d\mu_{mn}(x)=&\int_0^{+\infty} \PP \Bigl(x:\ -\ln q_m(x) \geq t \Bigr) \ dt \\=
&\int_{0}^{+\infty} \PP \Bigl(x:\ q_m(x) \leq e^{-t} \Bigr) \ dt \\ \leq \ &\int_{0}^{+\infty} \exp\left\{ {m \over 2} 
\left(1- e^{-t} -t \right) \right\} \ dt. \endsplit$$
Now,
$$\split \int_0^{+\infty} \exp\left\{ {m \over 2} \left(1-e^{-t} - t \right) \right\} \ dt =
&\int_0^1 \exp\left\{ {m \over 2} \left(1-e^{-t} - t \right) \right\} \ dt \\&\qquad+\int_1^{+\infty} 
\exp\left\{ {m \over 2} \left(1-e^{-t} - t \right) \right\} \ dt  \\ 
\leq \ &\int_0^1 \exp\left\{ - {m t^2 \over 6} \right\} \ dt + \int_0^{+\infty} \exp\left\{ -{m t \over 2} \right\} \ dt \\
\leq \ &\sqrt{3\pi \over 2m} + {2 \over m}.\endsplit$$
Summarizing, 
$$\EE \left| \ln q_m \right| \ \leq \ \sqrt{\pi \over m} + \sqrt{3 \pi \over 2m} + {2 \over m} \ < {6 \over \sqrt{m}}$$
and the proof of Part (3) follows.
{\hfill \hfill \hfill} \qed
\enddemo

\proclaim{(5.2) Theorem} Let $\alpha_1, \ldots, \alpha_k$ be non-negative reals such that 
$\alpha_1+ \ldots + \alpha_k=1$, let $Q_1, \ldots, Q_k$ be $n \times n$ positive definite matrices and let 
$m$ be a positive integer.
Then
$$\split &\max \Sb X \succeq 0  \\ \tr(X)=1 \\ \rk X \leq m \endSb \sum_{i=1}^k \alpha_i \ln \langle Q_i, X \rangle \ \leq \ 
\max \Sb X \succeq 0 \\ \tr(X)=1 \endSb \sum_{i=1}^k \alpha_i \ln \langle Q_i, X \rangle \\ &\qquad \qquad \leq \ {15 \over \sqrt{m}} + \max \Sb X \succeq 0 \\ 
\tr(X)=1 \\ \rk(X) \leq m  \endSb \sum_{i=1}^k \alpha_i \ln \langle Q_i, X \rangle. \endsplit$$
\endproclaim
\demo{Proof} The first inequality obviously holds.

Let $A$ be a matrix where the maximum value of the function 
$$X \longmapsto \sum_{i=1}^k \alpha_i \ln \langle Q_i, X \rangle$$
is attained on the set $X$ of positive semidefinite matrices of trace 1.
Rescaling $Q_i \longrightarrow \tau_i Q_i$ for some positive $\tau_1, \ldots, \tau_k$ if necessary, we may assume 
that $\langle Q_i , A \rangle =1$ for $i=1, \ldots, k$ and hence
$$\max \Sb X \succeq 0 \\ \tr(X)=1 \endSb \sum_{i=1}^k \alpha_i \ln \langle Q_i, X \rangle =
 \sum_{i=1}^k \alpha_i \ln \langle Q_i, A \rangle =0. \tag5.2.1$$
 Since $A$ is positive semidefinite, we can write $A=T^2$ for some symmetric
$n \times n$ matrix $T$. 

Let us fix the standard Gaussian probability measure $\mu_n$ in ${\Bbb R}^n$ with density 
$${1 \over (2\pi)^{n/2}} e^{-\|x\|^2/2}$$
and let $x_1, \ldots, x_m \in {\Bbb R}^n$ be $m$ independent random vectors.
Then 
$$\EE \| Tx_j\|^2 =\EE \langle Tx_j,\ Tx_j \rangle = \EE \langle T^2 x_j,\ x_j \rangle = \tr\left(T^2 \right)=
\tr(A)=1.$$
Applying Part (1) of Lemma 5.2, we conclude that 
$$\PP\left(x_1, \ldots, x_m:\ {1 \over m} \sum_{j=1}^m\|Tx_j\|^2 \ \geq \ 1+ {3 \over \sqrt{m}} \right)
 \ \leq \ \exp\left\{ -{9 \over 8}\right\} \ < \ 0.33
\tag5.2.2$$
(we use that $\ln (1+s) \leq s -s^2/4$ for $0 \leq s \leq 1$).

Let us define quadratic forms
$$q_i(x)=\langle Q_i, x \otimes x \rangle \quad \text{for} \quad i=1, \ldots, k.$$
Then
$$\split \EE q_i(Tx_j)= &\langle Q_i Tx_j,\ Tx_j \rangle =\langle T Q_i T x_j,\ x_j \rangle =\tr \left(T Q_i T \right)\\=
&\tr \left(Q_i T^2 \right) =\langle Q_i, A \rangle =1 \quad \text{for} \quad i=1, \ldots, k. \endsplit$$
Therefore, by Part (3) of Lemma 5.1, 
$$\EE \left| \ln \left({1 \over m} \sum_{j=1}^m q_i(Tx_j)\right)\right| \ \leq \ {6 \over \sqrt{m}} \quad \text{for} \quad i=1, \ldots, k$$ 
and hence 
$$\EE \left| \sum_{i=1}^k \alpha_i \ln \left({1 \over m} \sum_{j=1}^m q_i(Tx_j) \right) \right| \ \leq \ {6 \over \sqrt{m}}.$$
Therefore, by the Markov inequality,
$$\PP \left(x_1, \ldots, x_m:\quad \sum_{i=1}^k \alpha_i \ln \left({1 \over m} \sum_{j=1}^m q_i(Tx_j) \right)\ \leq \ 
-{12 \over \sqrt{m}}  \right) \ \leq \ 0.5. \tag5.2.3$$
From (5.2.2)--(5.2.3) we conclude that there are points $x_1, \ldots, x_m \in {\Bbb R}^n \setminus \{0\}$ such that 
$${1 \over m} \sum_{j=1}^m \|Tx_j\|^2  \ \leq \ 1+{3 \over \sqrt{m}}
\quad \text{and} \quad \sum_{i=1}^k \alpha_i \ln \left( {1 \over m} \sum_{j=1}^m q_i(Tx_j)\right) \ \geq \ -{12 \over \sqrt{m}}.$$
Let us define a matrix $Y$ by 
$$Y= \left(\sum_{j=1}^m \|Tx_j\|^2\right)^{-1} \sum_{j=1}^m \left(Tx_j\right) \otimes \left(Tx_j\right).$$
Then 
$$Y \succeq 0, \quad \tr(Y)=1, \quad \rk Y \leq m$$
and 
$$\split  \sum_{i=1}^k \alpha_i \ln \langle Q_i, Y \rangle =&\sum_{i=1}^k \alpha_i \ln \left({1 \over m} \sum_{j=1}^m q_i(x_j) \right) 
-\ln \left({1 \over m} \sum_{j=1}^m \|Tx_j\|^2 \right)\\
\ \geq\ -&{12 \over \sqrt{m}} - \ln \left(1 + {3 \over \sqrt{m}} \right) \ > \ -{15 \over \sqrt{m}}, \endsplit$$
and, in view of (5.2.1), the proof follows.
{\hfill \hfill \hfill} \qed
\enddemo

\subhead (5.3) Proof of Theorem 1.2 \endsubhead
 As in the proof of Theorem 1.1 in Section 4, without loss of generality we assume that
$$\sum_{i=1}^k Q_i =I. \tag5.3.1$$
Since $a \in \conv\left(\psi\left({\Bbb  R}^n\right)\right)$, we can write 
$$a_i =\langle Q_i, X \rangle \quad \text{for} \quad i=1, \ldots, k$$ and 
some $X \succeq 0$. Moreover, in view of (5.3.1), we have 
$$1=\sum_{i=1}^k a_i = \left\langle \sum_{i=1}^k Q_i,\ X \right\rangle = \langle I,\ X \rangle =\tr(X).$$
We note that 
$$\sum_{i=1}^k a_i \ln  \langle Q_i, X \rangle = \sum_{i=1}^k a_i \ln a_i.$$
By Theorem 5.2, there is a $n \times n$ symmetric matrix $Y$, such that $Y \succeq 0$, $\rk Y \leq m$ and 
$${15 \over \sqrt{m}} + \sum_{i=1}^k a_i \ln \langle Q_i, Y \rangle \ > \  \sum_{i=1}^k a_i \ln a_i.$$
Let
$$b_i = \langle Q_i, Y\rangle \quad \text{for} \quad i=1, \ldots, k.$$
Then
$$\sum_{i=1}^k b_i =   \left\langle \sum_{i=1}^k Q_i,\ Y\right\rangle = \tr(Y)=1.$$ 
Since $\rk Y \leq m$, we can write
$$Y={1 \over m} \sum_{j=1}^m y_j \otimes y_j$$
for some $y_1, \ldots, y_m \in {\Bbb R}^n$. Then 
$$b_i ={1 \over m} \sum_{j=1}^m q_i\left(y_j\right) \quad \text{for} \quad i=1, \ldots, k$$ and 
$b$ is a convex combination of at most $m$ points from $\psi\left({\Bbb R}^n\right)$.
{\hfill \hfill \hfill} \qed

\head Acknowledgment \endhead

I am grateful to Roman Vershynin for several useful suggestions and references.

\enddocument
\end